\documentclass[oneside,english,11pt]{amsart}
\usepackage{lmodern}
\usepackage[T1]{fontenc}
\usepackage[latin9]{inputenc}
\usepackage{geometry}
\usepackage{amsbsy}
\usepackage{amstext}
\usepackage{amsthm}
\usepackage{amssymb}
\usepackage{color}
\usepackage{nicefrac}
\usepackage[noend]{algorithmic}
\usepackage[ruled,vlined]{algorithm2e}
\usepackage{caption}

\makeatletter
\newcommand{\RemoveAlgoNumber}{\renewcommand{\fnum@algocf}{\AlCapSty{\AlCapFnt\algorithmcfname}}}
\newcommand{\RevertAlgoNumber}{\algocf@resetfnum}
\numberwithin{equation}{section}
\numberwithin{figure}{section}
\theoremstyle{plain}
\newtheorem{thm}{\protect\theoremname}
  \theoremstyle{plain}
  \newtheorem{prop}[thm]{\protect\propositionname}
  \theoremstyle{plain}
  \newtheorem{lem}[thm]{\protect\lemmaname}
  \theoremstyle{plain}
  \newtheorem{cor}[thm]{\protect\corollaryname}
  \theoremstyle{remark}
  \newtheorem{rem}[thm]{\protect\remarkname}

\DeclareMathOperator*{\argmin}{arg\,min}

\makeatother

\usepackage{babel}
  \providecommand{\corollaryname}{Corollary}
  \providecommand{\lemmaname}{Lemma}
  \providecommand{\propositionname}{Proposition}
  \providecommand{\remarkname}{Remark}
\providecommand{\theoremname}{Theorem}

\begin{document}
	\RemoveAlgoNumber
\global\long\def\e{e}
\global\long\def\bs{\boldsymbol{\sigma}}
\global\long\def\bx{\mathbf{x}}
\global\long\def\by{\mathbf{y}}
\global\long\def\bv{\mathbf{v}}
\global\long\def\bu{\mathbf{u}}
\global\long\def\bn{\mathbf{n}}
\global\long\def\grad{\nabla}
\global\long\def\Hess{\nabla^{2}}
\global\long\def\ddq{\frac{d}{dR}}
\global\long\def\qs{q_{\star}}
\global\long\def\qss{q_{\star\star}}
\global\long\def\Es{E_{\star}}
\global\long\def\EH{E_{\Hess}}
\global\long\def\Esh{\hat{E}_{\star}}
\global\long\def\ds{d_{\star}}
\global\long\def\Cs{\mathscr{C}_{\star}}
\global\long\def\nh{\boldsymbol{\hat{\mathbf{n}}}}
\global\long\def\BN{\mathbb{B}^{N}}
\global\long\def\ii{\mathbf{i}}
\global\long\def\SN{\mathbb{S}^{N-1}}
\global\long\def\SNq{\mathbb{S}^{N-1}(q)}
\global\long\def\SNqd{\mathbb{S}^{N-1}(q_{d})}
\global\long\def\SNqp{\mathbb{S}^{N-1}(q_{P})}
\global\long\def\nd{\nu^{(\delta)}}
\global\long\def\nz{\nu^{(0)}}
\global\long\def\cls{c_{LS}}
\global\long\def\qls{q_{LS}}
\global\long\def\dls{\delta_{LS}}
\global\long\def\E{\mathbb{E}}
\global\long\def\P{\mathbb{P}}
\global\long\def\R{\mathbb{R}}
\global\long\def\spp{{\rm Supp}(\mu_{P})}
\global\long\def\indic{\mathbf{1}}
\newcommand{\SNarg}[1]{\mathbb S^{N-1}(#1)} 

\title{Following the ground-states of full-RSB spherical spin glasses}

\author{Eliran Subag}
\begin{abstract}
We focus on spherical spin glasses whose Parisi distribution has support of the form
$[0,q]$. For such models we construct paths from the origin to the sphere which consistently remain close to the ground-state
energy on the sphere of corresponding radius. The construction uses
a greedy strategy, which always follows a direction corresponding to the most
negative eigenvalues of the Hessian of the Hamiltonian. For finite mixtures $\nu(x)$ it provides 
an algorithm of time complexity $O(N^{\deg(\nu)})$ to find w.h.p. points with the ground-state energy, up to a small error.

For the pure spherical models, the same algorithm reaches the energy
$-E_{\infty}$, the conjectural terminal energy for gradient descent.
Using the TAP formula for the free energy, for full-RSB models with
support $[0,q]$, we are able to prove the correct lower bound on
the free energy (namely, prove the lower bound from Parisi's formula),
assuming the correctness of the Parisi formula only in the replica
symmetric case. 
\end{abstract}

\maketitle

\section{Introduction}

Gaussian processes on the sphere in $\R^N$ are natural models for the study of smooth random functions in high dimensions, in general, and their (non-convex) optimization, in particular. The spherical spin glass models, which we consider in this work,  are  such processes, or Hamiltonians, defined on the sphere of radius $\sqrt N$ by  polynomials with independent Gaussian coefficients. As sequences of processes in the dimension, they are characterized by the property that their covariance function only depends on the angle between two points, up to scaling by $N$.\footnote{See \eqref{eq:cov} below for the explicit formula for the covariance. By \cite{Schoenberg}, this is the most general form of covariance that only depends on the angle (up to also allowing the $p=0,1$ terms in the definition of the mixture $\nu(x)$).}
 
The `order parameter' of spherical spin glasses is the minimizing measure from Parisi's formula for the free energy, and it generically coincides with the limiting annealed overlap distribution. Namely, the law of the normalized inner product of two independent samples from the Gibbs measure.  
We will (mainly, but not only) be interested in models exhibiting full Replica Symmetry Breaking (RSB), for which the support of the latter measure is of the form $[0,q]$, at any temperature.

We will focus on the computational problem of efficient optimization when the objective function is the highly non-convex landscape of a full-RSB spherical spin glass. More precisely, we will construct a polynomial time algorithm whose input are the disorder coefficients, i.e., the Gaussian coefficients of the Hamiltonian, and whose output is, with high probability (w.h.p.), a point on the sphere where the value of the Hamiltonian is approximately its global minimum. 

The optimization algorithm will arise naturally from new insights we gain into the geometric structure of the set of ground-state configurations, i.e., the approximate minimizers of the energy. 
In fact, we will consider an extended Hamiltonian, defined also in the interior of the sphere, and for each radius smaller than $\sqrt N$ we will examine the set of ground-state configurations on the corresponding sphere. 

Our main interest will be to understand how one can `follow' those ground-states as the radius increases. Our approach will exploit (sometimes as a tool in the proofs, but mostly only as inspiration) several 
properties of the ultrametric tree: it is supported on a continuous range of overlaps (in the full-RSB case); it branches in many orthogonal directions; and its vertices, being centers of heavy spherical bands, are approximate ground-state configurations on the sphere of corresponding radius. 
We will show that using a local greedy strategy, following a direction of the most negative eigenvectors of the Hessian on the orthogonal plane at each step, we are able to imitate the evolution of a path on the tree. In particular, the path constructed in this way, starting from the origin and ending on the sphere of radius $\sqrt N$, on which the original Hamiltonian is defined, consistently remains close to the ground-state energy of the corresponding radius.

\subsection{Model definition and background}

Let $J_{i_{1},...,i_{p}}^{(p)}$ be i.i.d. standard Gaussian variables,
and let $\gamma_{p}\geq0$, $p\geq2$, be a deterministic sequence
decaying exponentially in $p$. The \emph{Hamiltonian} of the spherical\emph{
}$p$-spin spin glass model with \emph{mixture} $\nu(x)=\sum_{p\geq2}\gamma_{p}^{2}x^{p}$
is
\begin{equation}
H_{N}\left(\bs\right)=\sum_{p=2}^{\infty}\frac{\gamma_{p}}{N^{\left(p-1\right)/2}}\sum_{i_{1},...,i_{p}=1}^{N}J_{i_{1},...,i_{p}}^{(p)}\sigma_{i_{1}}\cdots\sigma_{i_{p}},\quad\boldsymbol{\sigma}=(\sigma_{1},\ldots,\sigma_{N})\in\SN,\label{eq:Hamiltonian}
\end{equation}
where $\SN$ denotes the (Euclidean) sphere of radius $\sqrt{N}$
in dimension $N$. Equivalently, it is the centered Gaussian process
with covariance function 
\begin{equation}
	\label{eq:cov}\E \big\{
	H_{N}(\bs)H_{N}(\bs') \big\} = N\nu(R(\bs,\bs')),
\end{equation} where $R(\bs,\bs')=\bs\cdot\bs'/N$
is the usual overlap function. 

The associated Gibbs measure is the random probability measure on
$\SN$ given by
\[
\frac{dG_{N,\beta}}{d\bs}(\bs):=\frac{1}{Z_{N,\beta}}e^{-\beta H_{N}(\bs)},\quad Z_{N,\beta}:=\int_{\SN}e^{-\beta H_{N}(\bs)}d\bs,
\]
where $d\bs$ denotes the normalized Haar measure on the sphere. The
normalization factor $Z_{N,\beta}$ and $F_{N,\beta}=\frac{1}{N}\log Z_{N,\beta}$
are the usual partition function and free energy, respectively. 

The limiting free energy is given by the celebrated Parisi formula
\cite{ParisiFormula,Parisi1980,Crisanti1992,Talag,Chen},
\begin{equation}
F_{\beta}:=\lim_{N\to\infty}\E F_{N,\beta}=\inf_{x}\mathcal{P}(x),\label{eq:Parisi}
\end{equation}
where the infimum is over all distribution functions $x$ on $[0,1]$,
$\hat{x}(q)=\int_{q}^{1}x(s)ds$ and, with arbitrary $\hat{q}$ such
that $x(\hat{q})=1$, the Crisanti-Sommers representation \cite{Crisanti1992},
 of the functional
$\mathcal{P}(x)$ is given by 
\[
\mathcal{P}(x)=\frac{1}{2}\Big(\int_{0}^{1}x(q)\beta^{2}\nu'(q)dq+\int_{0}^{\hat{q}}\frac{dq}{\hat{x}(q)}+\log(1-\hat{q})\Big).
\]
The infimum in (\ref{eq:Parisi}) is attained uniquely at $x_{P}=x_{P,\beta}$,
called the Parisi distribution. Denote by $S_{P}=S_{P,\beta}$ the
support of the measure corresponding to $x_{P}$. 

As first shown by Auffinger and Chen \cite[Theorem 6]{AuffingerChen2015},
interestingly, in the full-RSB case the Parisi distribution in fact
simplifies. Namely, they showed that at any interior point of $S_{P}$
the density of $x_{P}$ takes a very specific form (see (\ref{eq:eta})
below). We collect in the following proposition several well-known
related results from the works of Chen and Panchenko \cite{PaCheChaos},
Jagannath and Tobasco \cite{JagannathTobascoBdsCplxSph} and Talagrand
\cite{Talag}. We will mainly be interested in models satisfying the (equivalent) conditions in the proposition.
\begin{prop}[full-RSB models \cite{PaCheChaos,JagannathTobascoBdsCplxSph,Talag}\footnote{For $\beta\leq\nu''(0)^{-\frac{1}{2}}$, Point (\ref{enu:FRSB3})
		follows from Point (\ref{enu:FRSB1}) by \cite[Proposition 2.3]{Talag}
		(see the proof of Theorem \ref{thm:FEwoutParisi} below). For $\beta>\nu''(0)^{-\frac{1}{2}}$
		Point (\ref{enu:FRSB3}) follows from Point (\ref{enu:FRSB1}) by
		\cite[Proposition 2]{PaCheChaos}. From (\ref{eq:Parisi}) and Point
		(\ref{enu:FRSB2}) one can see that $q_{P}\to1$ as $\beta\to\infty$.
		Thus, Point (\ref{enu:FRSB1}) follows from Point (\ref{enu:FRSB2})
		by \cite[Corollary  1.6]{JagannathTobascoBdsCplxSph}.}]
\label{prop:1} The following are equivalent:
\begin{enumerate}
\item \label{enu:FRSB1}$\nu''(q)^{-\frac{1}{2}}$ is concave on $(0,1]$.
\item \label{enu:FRSB2}For any $\beta>0$, $S_{P}$ is of the form $[0,q_{P}]$.
\item \label{enu:FRSB3}If $\beta\leq\nu''(0)^{-\frac{1}{2}}$, then $S_{P}=\{0\}$,
and if $\beta>\nu''(0)^{-\frac{1}{2}}$, then $S_{P}=[0,q_{P}]$ where
$q_{P}=q_{P,\beta}$ is the unique solution of $\nu''(q)^{-\frac{1}{2}}=\beta(1-q)$
and 
\begin{equation}
x_{P}=\begin{cases}
\frac{1}{\beta}\eta(q):=-\frac{1}{\beta}\frac{d}{dq}\nu''(q)^{-1/2}=\frac{\nu'''(q)}{2\beta\nu''(q)^{3/2}}, & \text{if }q\in[0,q_{P}),\\
1, & \text{if }q\in[q_{P},1].
\end{cases}\label{eq:eta}
\end{equation}
\end{enumerate}
\end{prop}

 In view of Point (\ref{enu:FRSB1}), we also mention that
by \cite[Corollary  1.6]{JagannathTobascoBdsCplxSph}, if one splits
$(0,1)$ into intervals on which $\nu''(q)^{-\frac{1}{2}}$ is concave
and intervals on which it is convex, then full-RSB can occur on the
former intervals and only there, and finite RSB can occur on the latter
intervals and only there. 

Obviously, by combining Parisi's formula with Point (\ref{enu:FRSB3})
one obtains the limiting free energy. Another way to express it, more
relevant to us, is by the generalized TAP formula recently proved in \cite{FElandscape}.
The latter states, for models as in Proposition \ref{prop:1}, that
for any $q\in[0,q_{P}]$,
\begin{equation}
F_{\beta}=\beta\Es(q)+\frac{1}{2}\log(1-q)+F_{\beta}(q),\label{eq:TAP}
\end{equation}
where, denoting by 
\[
\SNq=\left\{ \bs:\,\|\bs\|=\sqrt{Nq}\,\right\} 
\]
the sphere of radius $\sqrt{Nq}$, 
\begin{equation}
-\Es(q):=\lim_{N\to\infty}\frac{1}{N}\E\min_{\bs\in\SNq}H_{N}(\bs)\overset{\text{a.s.}}{=}\lim_{N\to\infty}\frac{1}{N}\min_{\bs\in\SNq}H_{N}(\bs),\label{eq:GS}
\end{equation}
is the corresponding ground-state energy, (where $H_{N}(\bs)$ is
defined through (\ref{eq:Hamiltonian}) in the interior of $\SN$
as well) where $F_{\beta}(q)$ is the limiting free energy of the
mixture
\[
\nu_{q}(s)=\nu(q+(1-q)s)-\nu(q)-(1-q)\nu'(q)s,
\]
and where $\frac{1}{2}\log(1-q)$ is an entropy term, equal to the
logarithmic growth rate of the volume of the subset of points  $\bs\in\SN$ such that $R(\bs-\bs_0,\bs_0)=0$
for some arbitrary $\bs_0\in\SNq$. Moreover,
for $q=q_{P}$, $\nu_{q}(s)$ is replica symmetric and $F_{\beta}(q)=\beta^{2}\nu_{q}(1)/2$,
so that the only non-trivial part in (\ref{eq:TAP}) is the ground-state
energy $\Es(q)$. 

For the full-RSB models, one can obtain $\Es(1)$ directly from the Parisi
formula by substituting $x_{P}$, and taking the $\beta\to\infty$
limit of $F_{\beta}/\beta$. By straightforward calculus, this leads to the following proposition.
The proposition can be also derived from the $\beta\to\infty$ analogue of the Parisi formula  of  Chen and Sen \cite{ChenSen}
and Jagannath and Tobasco \cite{JagannathTobascoLowTemp} 
for general spherical
models, and it was, in particular, stated in Proposition 2 of \cite{ChenSen}.
\begin{prop}[Ground state energy, Chen-Sen \cite{ChenSen}]
\label{prop:GSChenSen} If $\nu''(q)^{-\frac{1}{2}}$
is concave on $(0,1]$, 
\[
\Es(1)=\int_{0}^{1}\nu''(t)^{1/2}dt.
\]
\end{prop}

By scaling $H_{N}$ from $\SNq$ to the standard sphere $\SN=\SN(1)$,
noting that the corresponding scaled mixture $s\mapsto\nu(qs)$ satisfies
the concavity condition above, one also obtains that, in the setting
of Proposition \ref{prop:GSChenSen},
\begin{equation}
\Es(q)=\EH(q):=\int_{0}^{q}\nu''(t)^{1/2}dt.\label{eq:GSq}
\end{equation}

\subsection{\label{subsec:lowerbound}A greedy strategy based on $\protect\Hess H_{N}(\protect\bs)$: constructing paths of energy $-\EH(q)$ and
a general lower bound on $\protect\Es(q)$}

By the Borell-TIS inequality \cite{Borell,TIS}, w.h.p., for any $q\in[0,1]$ there exists some $\bs_{q}\in\SNq$ at which
$\frac{1}{N}H_{N}(\bs_{q})\approx-\Es(q)$. From \cite[Theorem 6]{FElandscape}, however, for models
as in Proposition \ref{prop:1} there is a whole \emph{continuous
path} in $q\in[0,1]$ satisfying the same. In fact, \emph{any path
in the infinite ultrametric tree} \cite{MPSTV2,MPSTV1,ultramet} has this property
and (when appropriately discretized) is, moreover, piece-wise given by short intervals where $\bs_{q+t}=\bs_{q}+\sqrt{Nt}\bv$,
for some $\bv$ unit vector orthogonal to $\bs_{q}$ (with infinitely
many such intervals, or levels to the tree, in the $N\to\infty$ limit). 

One may wish to gain insights from the above into how one can construct
paths that `follow' the ground-state energy $-\Es(q)$. Assuming the
path is of the form as above, this essentially boils down to the question: provided
that $\frac{1}{N}H_{N}(\bs_{q})\approx-\Es(q)$, how should we
choose a `good direction' $\bv$ in the orthogonal space  to  minimize the energy? 

For $\bs\neq0$, denoting by
\[
\nabla_{E}H_{N}(\bs)\text{ \,\,and \,\,}\nabla_{E}^{2}H_{N}(\bs)
\]
the usual Euclidean gradient and Hessian, define the projected gradient and Hessian 
\begin{equation}
\grad H_{N}(\bs)=M^{T}\nabla_{E}H_{N}(\bs)\text{ \,\,and \,\,}\Hess H_{N}(\bs)=M^{T}\nabla_{E}^{2}H_{N}(\bs)M,\label{eq:Hess}
\end{equation}
where $M=M^T=I-\bs\bs^T/\|\bs\|^2$ is the orthogonal projection matrix to the orthogonal space
of $\bs$. Since in our construction of the path $H_{N}(\bs_{q})$ should approximately be the global
minimum over $\SNq$, we should have $\grad H_{N}(\bs_{q})\approx0$, so we probably
should not choose $\bv$ based on the gradient. If this is the case, a natural way to proceed is to use the second order approximation and take $\bv$ such that
$\frac{1}{2}\bv^{T}\Hess H_{N}(\bs_{q})\bv$ is as negative as possible. 

For deterministic $\bs\in\SNq$, $\Hess H_{N}(\bs)$ is invariant
w.r.t. rotations of the orthogonal space to $\bs$, and the eigenvalues
of $\Hess H_{N}(\bs)$ are distributed as those of (see (\ref{eq:H2derivative}))
\[
\sqrt{(N-1)/N}\nu''(q)^{1/2}\mathbf{G},
\]
where $\mathbf{G}$ is a matrix of dimension $N-1$ from the Gaussian orthogonal ensemble
(GOE),\footnote{In this paper, a GOE matrix of dimension $N$ is a random matrix with independent centered,
	Gaussian entries up to symmetry, with variance $2/N$ on the diagonal
	and $1/N$ off the diagonal. Namely, the scaling we use is such that the
	limiting empirical measure is supported on $[-2,2]$.} plus an extra eigenvalue $\lambda=0$ that corresponds to the eigenvector $\bs$. In particular, typically,
\[
\min_{\|\bv\|=1}\frac{1}{2}\bv^{T}\Hess H_{N}(\bs)\bv\approx-\nu''(q)^{1/2},
\]
which is exactly what one would expect, in order to obtain (\ref{eq:GSq}).
The following lemma, which we will prove in Section \ref{sec:Proofs}, similarly controls the `edge'  of the spectrum of $\Hess H_N(\bs)$, 
uniformly over the ball of radius $\sqrt{N}$, instead of one deterministic point. Let $\lambda_{i}(\bs)$, $i=1,\ldots,N-1$,
be the eigenvalues of $\Hess H_{N}(\bs)$, with the $0$ eigenvalue that corresponds to the eigenvector $\bs$ removed, ordered in non-decreasing
order.
\begin{lem}[Uniform control on the edge]
\label{lem:good_eigvals} Let $\nu(x)$ be a general mixture. For any
$\epsilon>0$, there exist $\delta,\,c,\,K>0$ such that 
\begin{equation}
\label{eq:eigvalsbd}
\begin{aligned}
\P\Big\{\,\forall\bs \in \BN:\: \,  &\#\big\{i: \lambda_i(\bs) \leq -2\nu''(q)^{\frac{1}{2}}+\epsilon \big\} \geq N \delta,\\
 &\#\big\{i: \lambda_i(\bs) \leq -2\nu''(q)^{\frac{1}{2}}-\epsilon \big\} \leq  K \,\Big\}\geq1-e^{-Nc}	,
\end{aligned}
\end{equation}
where $\BN:= \big\{\bs : \|\bs\|^{2}/N\in(0,1] \big\}$ and $q=q(\bs):=\|\bs\|^{2}/N$.
\end{lem}


With the lemma at our disposal, we are ready to
construct paths that, in the full-RSB case, follow the evolution of
the ground-state energy. 
 Set $\bs_{0}=0$ and suppose that $\bv_{0},\ldots,\bv_{k-1}\in\R^{N}$
with $\|\bv_{i}\|=1$ are $k\geq1$ directions satisfying $\bs_{j/k}\cdot\bv_{j}=0$,
where we define $\bs_{j/k}=\sqrt{N/k}\sum_{i=0}^{j-1}\bv_{i}$. Extend
this sequence to a continuous path parameterized by $[0,1]$ by interpolating,
for any $t\in[0,1/k]$,
\begin{equation}
\bs_{j/k+t}=\bs_{j/k}+\sqrt{Nt}\bv_{j}\in\SN\big(\hspace{0.05em}j/k+t\big).\label{eq:interp}
\end{equation}
In Section \ref{sec:Proofs} we will prove the following theorem. The corresponding pseudo-code can also be found there.

\begin{thm}[Optimization algorithm]
	\label{thm:LB} Let $\nu(x)$ be a general mixture and let $\epsilon>0$
	and $k\geq1$. Then there exist constants $C$, $C'$, $\eta$, depending only on $\nu$, $\epsilon$ and $k$, such that $\lim_{\substack{\epsilon\to0\\k\to\infty}}\eta =0$ and:
	\begin{enumerate}
		\item\label{Item1} With probability at least $1-e^{-NC}$, for any choice of directions $\{\bv_j\}_{j=0}^{k-1}$ as above that,  for $\bv = \bv_j$, $\bs =\bs_{j/k}$, $q=\|\bs_{j/k}\|^2/N$ and any $j$, satisfy
		\begin{equation}
		\bv\cdot\grad H_{N}(\bs)\leq0\text{\quad\, and \quad\,}
		\bv^{T}\Hess H_{N}(\bs)\bv\leq -2\nu''(q)^{\frac{1}{2}}+\epsilon,\label{eq:vj}
		\end{equation}
		we have that 
		\[
		\forall q\in[0,1]:\ \ \frac{1}{N}H_{N}(\bs_{q})\leq-\EH(q)+\eta.
		\]
		\item \label{Item2} Assuming $\deg(\nu)<\infty$,\footnote{\label{ft:degnu}If $\deg(\nu)=\infty$, there are infinitely many disorder coefficients and just storing them cannot be done in finite complexity. In this case, one can either work with a finite mixture approximation by using the formula of \eqref{eq:Hamiltonian} with summation over $p$ only up to some large $p_0$ and uniformly bound the error w.h.p., or alternatively use queries of $H_N(\bs)$ to measure the complexity (in which case it is possible to show that $C'N$ queries are enough).} there exists an algorithm with  complexity $C' N^{\deg(\nu)}$ (measured in basic mathematical operations) that takes as input the disorder coefficients $J_{i_{1},...,i_{p}}^{(p)}$ and outputs directions $\bv_j$ as in  \eqref{Item1}, with probability at least $1-e^{-NC}$.
	\end{enumerate}
	\end{thm}
By combining Part \eqref{Item1} of Theorem \ref{thm:LB} and Lemma \ref{lem:good_eigvals} we immediately conclude the following. See Section \ref{sec:Proofs} for a short proof.
\begin{cor}
\label{cor:GSlb}Let $\nu(x)$ be a general mixture. Then, for any
$q$,
\begin{equation}
\label{eq:LB}
\Es(q)\geq\EH(q).
\end{equation}
\end{cor}

When $\nu''(q)^{-1/2}$ is concave, the lower bound of Corollary \ref{cor:GSlb} coincides with
the ground-state energy as follows directly from Parisi's formula or \cite{ChenSen,JagannathTobascoLowTemp},
see Proposition \ref{prop:GSChenSen} above. In contrast to
the latter, our proof of Corollary \ref{cor:GSlb} essentially
only uses well-known properties of GOE matrices and Gaussian fields, and in particular
does not rely on any information from Parisi's formula.

Note that on the event from \eqref{eq:eigvalsbd}, for any $\bs$ the set of vectors $\bv$ such that $\pm\bv/\|\bv\|$ satisfies the second inequality of \eqref{eq:vj} contains a linear subspace of dimension which is a positive fraction $N$ (omitting $\bv=0$). One of $\pm\bv/\|\bv\|$ also satisfies the first inequality of \eqref{eq:vj}. Hence, the directions $\bv_j$ as in Theorem \ref{thm:LB} can be chosen one after the other, starting with arbitrary $\bv_0$ and given  $\bs_{j/k}$, choosing $\bv_j$ so that the \eqref{eq:vj} holds and setting $\bs_{j+1/k}=\bs_{j/k}+\sqrt{N/k}\bv_j$.

\begin{rem}
In fact, by choosing many directions at each step instead of just one, one may also easily modify this construction and obtain a tree instead of a path, such that each branch has `direction' $\bv$ satisfying \eqref{eq:vj}. From the previous paragraph, the directions of the branches can be chosen to be orthogonal to each other, also across all levels of the tree. Moreover, the degree of each vertex can be taken to be diverging with $N$. By a diagonalization argument we can also make the number of levels  in the tree $k$ diverge, and the error $\eta$ vanish, as $N\to\infty$. In the full-RSB case, the leaves of this tree are near-ground-state configurations and the overlap of any two leaves is determined by their distance on the tree. In the large $N$ limit, the range of overlaps between all pairs of leaves covers the full range $[0,1]$. 
The fact that there exist near-ground-state configurations with approximately any overlap in $[0,1]$ was concluded by Auffinger and Chen \cite[Theorem 4]{AuffingerChenEnergyLandscape} for full-RSB models from a general principle they prove for spherical models (with even interactions), using the Parisi formula.
\end{rem}

\subsection{Using many orthogonal directions: the matching upper bound on $\protect\Es(q)$
in the full-RSB case}

In Section \ref{subsec:lowerbound} we discussed the lower bound on $\Es(q)$.
At the moment, when $\nu''(q)^{-1/2}$ is concave, the matching upper
bound to \eqref{eq:LB} is known to us directly from the Parisi formula, or  \cite{ChenSen,JagannathTobascoLowTemp} which build
on the formula. In this section we explain how the same bound
can be derived by a different approach, using the infinitary nature
(or, duplication property) of the ultrametric structure proved by
Panchenko in the seminal work \cite{ultramet} and using ideas from
\cite{FElandscape}. The basic result we deduce from \cite{FElandscape} is
the following lemma, which we prove in Section \ref{sec:Proofs}.
\begin{lem}
\label{lem:orthoGS}Let $\nu(x)$ be a general mixture, let $\beta>0$ be arbitrary, and assume
that $q,\,q'\in S_{P}\cap(0,1)$ and $q<q'$. Then, for some $\delta_{N}$
and $k_{N}$ tending to $0$ and $\infty$, respectively, with probability
tending to $1$ as $N\to\infty$, there exist points
\[
\bs\in\SNq\text{\,\,\,\,and\,\,\,\,}\bs_{i}\in\SN(q'),\:i=1,\ldots,k_{N},
\]
such that for any different $i,\,j\leq k_{N}$,
\begin{equation}
\big|R(\bs,\bs_{i}-\bs)\big|<\delta_{N}\text{\,\,\,\,and\,\,\,\,}\,\big|R(\bs_{i}-\bs,\bs_{j}-\bs)\big|<\delta_{N},\label{eq:OlapBd}
\end{equation}
and
\begin{equation}
\Big|\frac{1}{N}H_{N}(\bs)+\Es(q)\Big|<\delta_{N}\text{\,\,\,\,and\,\,\,\,}\Big|\frac{1}{N}H_{N}(\bs_{i})+\Es(q')\Big|<\delta_{N}.\label{eq:ValBd}
\end{equation}
\end{lem}

Assuming the \emph{approximate} \emph{orthogonality} (\ref{eq:OlapBd}),
there must be a large set of directions
\[
\bu_{i}:=\frac{\bs_{i}-\bs}{\|\bs_{i}-\bs\|}\approx\frac{\bs_{i}-\bs}{\sqrt{N(q'-q)}}
\]
out of the $k_{N}\gg1$ directions such that the projection of the
gradient $\grad H_{N}(\bs)$ onto each of them is small (from Pythagoras'
theorem). Therefore, if $q'-q$ is very small, for those directions,
from (\ref{eq:ValBd}) we have that 
\[
\Es(q')-\Es(q)\approx\frac{1}{N}\big(H_{N}(\bs)-H_{N}(\bs_{i})\big)\approx-\frac{q'-q}{2}\bu_{i}^{T}\Hess H_{N}(\bs)\bu_{i}.
\]
Hence, if we are able to let $q'\to q$ while remaining in $S_{P}$,
the (one-sided) derivative of $\Es(q)$ at $q$ can be larger than
$\nu''(q)^{1/2}$, only if w.h.p. $\Hess H_{N}(\bs)$ has many eigenvalues
smaller than $-2\nu''(q)^{1/2}$. This would contradict Lemma \ref{lem:good_eigvals}, from which we will conclude
the following proposition. See Section \ref{sec:Proofs} for the proof. We reiterate that here our argument relies
on a heavy tool: the ultrametricity property \cite{ultramet}. However,
it is of geometric nature, in contrast to the analytical nature of  the proof of Parisi's formula or \cite{ChenSen,JagannathTobascoLowTemp}.
\begin{prop}
\label{prop:upper bound}Let $\nu(x)$ be a general mixture and $\beta>0$ be arbitrary. If $[q,q+t)\subset S_{P}$
for some small $t>0$, then 
\begin{equation}
\frac{d}{dq}^{+}\Es(q):=\lim_{\epsilon\searrow0}\frac{\Es(q+\epsilon)-\Es(q)}{\epsilon}=\nu''(q)^{1/2}.\label{eq:rightder}
\end{equation}
In particular, if $[0,q]\subset S_{P}$, then (see (\ref{eq:GSq}))
\[
\Es(q)=\EH(q).
\]
\end{prop}

\subsection{The lower bound on $F_{\beta}$ in the full-RSB case}

If $\nu''(q)^{-1/2}$ is concave,
the correct lower bound for $F_{\beta}$ can be derived provided we
only assume the correct lower bound for such models in the replica
symmetric regime\,---\,namely, we assume the following proposition.
\begin{prop}[Talagrand \cite{Talag}]
\label{prop:RSFE} Assume that $\nu''(x)^{-1/2}$
is concave on $(0,1]$. If 
\begin{equation}
\forall s\in(0,1):\,\,\beta^{2}\nu(s)+\log(1-s)+s<0,\label{eq:RSsolution}
\end{equation}
then $F_{\beta}=\frac{1}{2}\beta^{2}\nu(1)$.
\end{prop}

In \cite[Proposition 2.1]{Talag}, Talagrand proved a characterization
for the optimizer in Parisi's formula. Using it he showed in \cite[Proposition 2.3]{Talag}
that for general $\nu(x)$, $F_{\beta}=\frac{1}{2}\beta^{2}\nu(1)$
if and only if (\ref{eq:RSsolution}) holds.\footnote{To be precise, $F_{\beta}=\frac{1}{2}\beta^{2}\nu(1)$
if and only if  (\ref{eq:RSsolution}) holds as  non-strict inequality. Talagrand
did not show that if we have equality in (\ref{eq:RSsolution}) then
$F_{\beta}=\frac{1}{2}\beta^{2}\nu(1)$. However, this can be easily
verified by working with $\beta'<\beta$, so that 
 (\ref{eq:RSsolution}) holds as a strict inequality, and then letting $\beta'\to\beta$.} Above we stated only the weaker version we need. Nevertheless, one
should note that even the proof of this weaker version is based on
the general Parisi formula.

Combining Proposition \ref{prop:RSFE} with Corollary \ref{cor:GSlb}
and the TAP formula (\ref{eq:TAP}) lower bound obtained in \cite{FElandscape}\,---\,the proof of which does not rely neither on the Parisi
formula nor ultrametricity\footnote{This fact relies on a forthcoming improvement to \cite{FElandscape}.}\,---\,we will prove the following theorem in Section \ref{sec:Proofs}. We remark that
its proof, of course, does not rely on the matching upper bound of
Proposition \ref{prop:upper bound}.
\begin{thm}
\label{thm:FEwoutParisi}Assume that $\nu''(q)^{-\frac{1}{2}}$ is
concave on $(0,1]$. Then, provided that Proposition \ref{prop:RSFE}
holds, with $q_P$ and $x_P$ as in Proposition \ref{prop:1}, for any $\beta>0$,
\begin{equation}
F_{\beta}\geq\beta\int_{0}^{q_{P}}\nu''(t)^{1/2}dt+\frac{1}{2}\log(1-q_{P})+\frac{1}{2}\beta^{2}\nu_{q_{P}}(1).\label{eq:FElb}
\end{equation}
Moreover, since the right-hand side of (\ref{eq:FElb}) coincides
with $\mathcal{P}(x_{P})$, we recover the Parisi lower bound
\[
F_{\beta}\geq\inf_{x}\mathcal{P}(x).
\]
\end{thm}

\section{\label{sec:Algorithms}Remarks on  related models}
 
In this section we remark about several models related to the full-RSB spherical models. In particular, we discuss the optimization algorithm in the context of the 1-RSB pure spherical models.  

\subsubsection*{Ising spins}
Models with Ising spins are defined similarly to the spherical ones, with parameter space $\Sigma_N=\{\pm1\}^N$ instead of the sphere.
Very shortly after the first version of the current paper was posted, Montanari \cite{AndreaSK} proved that the  Sherrington-Kirkpatrick model (i.e., $\xi(x)=x^2$ with Ising spins) can be optimized in time complexity $O(N^2)$, conditionally on the assumption that the corresponding Parisi measure is full-RSB (which is widely believed to be the case). Montanari  used a message passing algorithm with orthogonal updates --- similarly to the orthogonality condition we impose $\bs_{j/k}\cdot\bv_{j}=0$ on the directions of the path \eqref{eq:interp}). Moreover, he showed that the algorithm can be used to construct approximate solutions to the TAP equations.

For the full-RSB spherical models we used the fact that  for any $q\in[0,1]$, points of the ultrametric tree are
the centers of `heavy' spherical bands (w.r.t. the Gibbs measure), or TAP states, and by appealing
to results from \cite{FElandscape}, are therefore ground-state configurations on $\SNq$. The analogue for full-RSB models with Ising spins
are
the generalized TAP states, introduced and computed by Chen, Panchenko
and the author in \cite{TAPChenPanchenkoSubag,TAP2ChenPanchenkoSubag}, which
as $q$ approaches $1$, approach the ground-state
configurations on the hyper-cube. We expect that efficient optimization for full-RSB models with Ising spins can be achieved by constructing paths from the origin to $\Sigma_N$ that maximize the energy plus the generalized TAP correction as $q$ increase, using orthogonal increments.
This will be investigated in future
work.

\subsubsection*{The CREM} In connection with the role played by concavity of $\nu''(q)^{-1/2}$
 in our analysis, we mention the following analogy with the Continuous
 Random Energy Model (CREM)\,---\,a random potential defined
 on the binary tree of depth $N$ with continuous variance profile
 (see \cite{BovKurk1,BovKurk2}). Recently, Addario-Berry and Maillard
 \cite{ABM_CREMalg} studied a greedy algorithm to find low energy
 states for the CREM which, starting at the root, jumps at each step
 to the descendant of minimal energy among all descendants at a given
 small depth relative to the current position. Their main result proves that the energy achieved by this algorithm is the
 algorithmic hardness energy threshold, i.e., the deepest energy that
 can be reached in polynomial time in an appropriate sense. Moreover,
 as they observed, this threshold coincides with the ground-state energy
 if and only if the variance profile of the CREM is a concave function.
 Lastly, we mention that the limiting overlap distribution of the CREM
 is of the form $[0,q]$ at any temperature if and only if the variance
 profile is strictly concave (see \cite[Theorem 3.6]{BovKurk2}).

\subsubsection*{Locally following ground-states}
 Our analysis uses in a crucial way the fact that for spherical full-RSB
 models we can locally `follow' the ground-state configurations as
 the radius $\sqrt q$ increases. Excluding the pure models for which this can
 be done trivially due to homogeneity, it seems to us that the same should not
be possible for mixed spherical models with finite RSB. For example,
 for some 1-RSB models it was shown in \cite[Corollary 11]{geometryMixed}
 that ground-state configurations of $\SN(q)$ are orthogonal to those
 of $\SN(q+\epsilon)$, for large $q$ and small $\epsilon$. Moreover, in the finite
 RSB case, overlap gaps for the set of near-ground-state configurations
 on the original sphere $\SN$ were proved in \cite{AuffingerChenEnergyLandscape}.
 Those phenomena are also intimately related to temperature chaos \cite{PaCheChaos,PaChaos},
 see specifically \cite[Theorem 3]{geometryGibbs}, \cite[Theorem 4]{geometryMixed}
 and \cite[Corollary 15]{FElandscape} whose proofs rely on related
 ideas.

We mention that one may also consider the problem of following the ground-states of stationary Gaussian fields in $\R^{N}$ in a
potential well. In fact, the replica computations of Fyodorov
 and Sommers \cite{FYODOROV2007128} (see also \cite{FyodorovLeDoussalHessian})
 suggest that for certain `long range' covariances the Parisi distribution takes a similar
 form to (\ref{eq:eta}), and is supported on an interval $[q_0,q_{EA}]$ that does not contain $0$. 
 As for the extended Hamiltonian we considered in the spherical case, for any radius $\sqrt q$ the restriction of the stationary Hamiltonian to the corresponding sphere is a spherical model with some $q$ dependent mixture (that can be explicitly calculated), but which in the stationary case will have a $1$-spin (or external field) term. 
 
\subsubsection*{The pure spherical models}

The pure spherical models $\nu(x)=x^{p}$ are special in several
respects and it will be instructive to inspect our optimization algorithm in their setting more carefully. First
and foremost, interestingly, for pure models $\EH(1)$ coincides with
a known and meaningful quantity
\[
E_{\infty}:=2\sqrt{(p-1)/p}=\sqrt{p(p-1)}\int_{0}^{1}s^{\frac{p-2}{2}}ds=\EH(1).
\]

The critical points of a pure model $H_{N}(\bs)$ on $\SN$ are
rather well-understood by the computation of means of Auffinger, Ben
Arous and {\v{C}}ern{\'y} \cite{A-BA-C}, an application of the second
moment method \cite{2nd}, and the convergence of the extremal process
proved by Zeitouni and the author \cite{pspinext} (also see \cite{ABA2,geometryMixed}
for mixed models). In particular, it is known that the (normalized) energy $-E_{\infty}$
is the threshold below which all critical points are of finite index,
and above which of diverging index. Since $H_{N}(\bs)$ is a.s. Morse, for almost every initial point on $\SN$, negative gradient flow converges to a local minimum, (a critical point of index
$0$) and its terminal energy has to be $-E_{\infty}$ or less.
Moreover, a randomly chosen local minimum of $H_{N}(\bs)$ on $\SN$
will have energy $-E_{\infty}+o(1)$ with overwhelming probability.
For this reason, $-E_{\infty}$ has been conjectured to be the energy
at which gradient descent typically terminates. 

Analyzing the time it takes for the flow to reach the minimum point, however, is a difficult problem which requires one to control regions of shallow slope along
the path of descent. Somewhat related aspects of Langevin dynamics,
a noisy version of gradient flow, like the mixing time, the aging
phenomenon and the Cugliandolo-Kurchan equations have been studied
in \cite{BDG,BDG1,BenArousGheissariJagannathBoundingflows,BenArousJagannathSpectralgap,Crisanti1993,PhysRevLett.71.173,Dembo2007,GheissariJAgannathSpectralgap}
(also see \cite{BenArousAging,BCKMoutofequilibrium,CugliandoloCourse,GuionnetDynamics}
for general surveys on this rich topic). However, to the best of our
knowledge, no rigorous upper bound is known for the typical time to
get close to $-E_{\infty}$ using negative gradient flow. 

Another property special to the pure case is that $H_{N}(\bs)$ is
a homogeneous function in $\R^N$, which allows one to run the optimization algorithm directly
on $\SN$. 
Precisely, we may start from an arbitrary point $\bs_0\in\SN$, and at each step find a direction $\bv$ as in the algorithm above and set\[
\bs_{i+1} = \bs/\|\bs\|, \mbox{\ \ \  where\ \ \ } \bs=\bs_i+\sqrt {N \tau} \bv,
\]
for some small parameter $\tau>0$.\footnote{If $\tau$ is normalized appropriately at each step, the path we obtain is simply the projection to $\SN$ of the path obtained from the algorithm of Theorem \ref{thm:LB} (see the proof of the theorem for a description of the algorithm).}
For this sequence on $\SN$,
one can easily verify that w.h.p.
\[
\frac{H_{N}(\bs_{k})}{N}\leq \frac{1}{(1+\tau)^{\frac{pk}{2}}}\frac{H_{N}(\bs_{1})}{N}-\nu''(1)^{1/2}\tau\frac{1-(1+\tau)^{-\frac{pk}{2}}}{(1+\tau)^{\frac{p}{2}}-1}+\tilde\eta,
\]
where $\tilde\eta$ is a constant that goes to $0$ as $\epsilon\to0$ (and depends on $\nu$ through bounds on the third order derivatives as in Lemma \ref{lem:Cont} below).
In particular, as we increase the number
of steps $k$, the initial value $H_{N}(\bs_{1})$ is washed away.
And since
\[
\lim_{\tau\to0}\frac{\tau}{(1+\tau)^{\frac{p}{2}}-1}=\frac{2}{p},
\]
in the large $k$ and small $\tau$ and $\epsilon$ limit, we obtain
\[
\frac{H_{N}(\bs_{k})}{N}\lesssim-\frac{2}{p}\nu''(1)^{1/2}=-E_{\infty}.
\]


\section{\label{sec:Proofs}Proofs}

We will need the following lemma in the proofs, which we take from \cite[Corollary 59]{geometryMixed}.
\begin{lem}[Order of derivatives \cite{geometryMixed}]
	\label{lem:Cont}  Let $\nu(x)$ be a general mixture. For appropriate
	$R,\,c>0$,  for $i=1,2,3$,
	\begin{align}
	\label{eq:Cont}
	\P \Big\{ \forall \bs\in\BN,\,\forall \bv \mbox{\ s.t.\ } \|\bv\|= 1: 	
	\big|\partial_{\bv}^{i}H_{N}(\bs)\big|<RN^{1-\frac{i}{2}} \Big\}&\geq 1-e^{-Nc},\\
	\P \Big\{ \forall \bs,\bs'\in\BN : 	
	\|\nabla_{E}^{2}H_{N}(\bs)-\nabla_{E}^{2}H_{N}(\bs')\|_{op}<\frac{R}{\sqrt{N}} \|\bs-\bs'\|     \Big\}&\geq 1-e^{-Nc},\label{eq:HessLip}
	\end{align}
	where $\partial_{\bv}^i$ is the  $i$-th order directional derivative in direction $\bv\in\R$ and $\|\cdot\|_{op}$ is the operator norm. 
\end{lem}

\subsection{Proof of Lemma \ref{lem:good_eigvals}}
By the min-max theorem, 
\[
\max_i|\hat{\lambda}_{i}(\bs)-\hat{\lambda}_{i}(\bs')|\leq\|\nabla_{E}^{2}H_{N}(\bs)-\nabla_{E}^{2}H_{N}(\bs')\|_{op},
\]
where
we denote by $\hat{\lambda}_{i}(\bs)$, $i=1,\ldots,N$,
the eigenvalues of  $\nabla_{E}^{2}H_{N}(\bs)$, ordered in non-decreasing
order.

Since $\nu''(q)^{1/2}$ is uniformly continuous
on $[0,1]$, we conclude that for small $t=t(\epsilon,\nu)>0$, if $T_{N}$
is an arbitrary $\sqrt{N}t$-net of $\BN$, then on the  event in \eqref{eq:HessLip} if
\begin{align}
\forall \bs\in T_N:&\ \  \#\big\{i: \hat\lambda_i(\bs) \leq -2\nu''(q)^{\frac{1}{2}}+\epsilon/2 \big\} \geq N \delta,\label{eq:TN1}
\\
\forall \bs\in T_N:&\ \ \#\big\{i: \hat\lambda_i(\bs) \leq -2\nu''(q)^{\frac{1}{2}}-2\epsilon \big\} \leq K,\label{eq:TN2}
\end{align}
then also 
\begin{align}
\forall \bs\in \BN:&\ \  \#\big\{i: \hat\lambda_i(\bs) \leq -2\nu''(q)^{\frac{1}{2}}+\epsilon \big\} \geq N \delta,\label{eq:BN1}
\\
\forall \bs\in \BN:&\ \ \#\big\{i: \hat\lambda_i(\bs) \leq -2\nu''(q)^{\frac{1}{2}}-\epsilon \big\} \leq K.\label{eq:BN2}
\end{align}

Recall that we defined $\Hess H_{N}(\bs)=M^{T}\nabla_{E}^{2}H_{N}(\bs)M$,
where  $M=I-\bs\bs^{T}/\|\bs\|^{2}$ is the orthogonal projection
matrix to the orthogonal space of $\bs$, and that $\lambda_{i}(\bs)$, $i=1,\ldots,N-1$,
are the eigenvalues of $\Hess H_{N}(\bs)$, with the $0$ eigenvalue that corresponds to the eigenvector $\bs$ removed, ordered in non-decreasing
order. By the mix-max theorem,\footnote{These inequalities sometime go by the name Cauchy's interlacing theorem, or Poincar\'{e}'s separation theorem.}
\begin{equation}
\hat{\lambda}_{i}(\bs)\leq\lambda_{i}(\bs)\leq\hat{\lambda}_{i+1}(\bs).\label{eq:hatlambda}
\end{equation}

Therefore, for large $N$, if \eqref{eq:TN1} and \eqref{eq:TN2} hold with $\hat{\lambda}_i(\bs)$ replaced by ${\lambda}_i(\bs)$, then they also hold in their original form with $\hat{\lambda}_i(\bs)$, where we may need to increase $K$. Similarly, for large $N$, if \eqref{eq:BN1} and \eqref{eq:BN2} hold, then they also hold with $\hat{\lambda}_i(\bs)$ replaced by ${\lambda}_i(\bs)$, where we may need to decrease $\delta$.

Combining the above we conclude that in order to complete the proof, it will be sufficient to show that for some $c>0$, for large $N$, with probability at least $1-e^{-cN}$,  \eqref{eq:TN1} and \eqref{eq:TN2} hold with $\hat{\lambda}_i(\bs)$ replaced by ${\lambda}_i(\bs)$. Of course, there exists a net $T_{N}$ as above with at most $e^{\rho N}$
elements, for large enough $\rho>0$.\footnote{E.g., by the fact that the minimal
number of balls of radius $\sqrt{N}t$ required to cover $\BN$  is bounded by the maximal number of disjoint balls
of radius $\sqrt{N}t/2$ with centers in $\BN$.}
Hence, by a union bound, to complete the proof it will be enough to show that for any $\hat c>0$, if $K$ is large enough and $\delta$ is small enough, then, for large $N$,
	\begin{align}
	&\sup_{\bs \in \BN}\P\Big\{ \#\big\{i: \lambda_i(\bs) \leq -2\nu''(q)^{\frac{1}{2}}+\epsilon \big\} < N \delta \Big\} \leq e^{-N\hat{c}},\label{eq:unifout1}
	\\
	&\sup_{\bs \in \BN}\P\Big\{ \#\big\{i: \lambda_i(\bs) \leq -2\nu''(q)^{\frac{1}{2}}-\epsilon \big\} >K\Big\} \leq e^{-N\hat{c}}.\label{eq:unifout2}
	\end{align}

As mentioned in the Introduction, for deterministic $\bs\in\SNq$, the joint law of
 $$\tilde\lambda_i(\bs):=\sqrt{N/(N-1)}\nu''(q)^{-1/2}\lambda_i(\bs)$$ is identical to the law of the eigenvalues of
a GOE matrix of dimension $N-1$. Since the law of $H_{N}(\bs)$
is invariant under rotations of $\R^{N}$, this can be seen by a direct
computation of the second order derivatives at the `north-pole' $\bn=(0,\ldots,0,\sqrt{Nq})$
using the formula for the Hamiltonian (\ref{eq:Hamiltonian}), 
\begin{equation}
\forall i,j<N:\,\,\frac{d}{dx_{i}}\frac{d}{dx_{j}}H_{N}\left(\bn\right)=\sum_{p=2}^{\infty}\frac{\gamma_{p}(Nq)^{\frac{p-2}{2}}}{N^{\frac{p-1}{2}}}\sum_{\{i_{1},...,i_{p}\}=\{i,j,N,\ldots,N\}}J_{i_{1},...,i_{p}}^{(p)}(1+\delta_{ij}),\label{eq:H2derivative}
\end{equation}
where the equality in the second summation is in the sense of multi-sets.

Denote by $\mu_{{\rm sc}}$ the semi-circle law, whose cumulative
distribution function is given by
\[
F_{{\rm sc}}(t)=\frac{1}{2\pi}\int_{-2}^{t}\indic_{|x|\leq2}\sqrt{4-x^{2}}dx,
\]
and choose some
\begin{equation}
\delta<F_{{\rm sc}}\Big(-2+\epsilon/\nu''(1)^{\frac{1}{2}}\Big)\leq F_{{\rm sc}}\Big(-2+\epsilon/\nu''(q)^{\frac{1}{2}}\Big)\,,\label{eq:delta}
\end{equation}
where the second inequality holds for $q\in(0,1]$. 

From our choice (\ref{eq:delta}), there exists a neighborhood $A$
of $\mu_{{\rm sc}}$ in the space of probability measures on $\R$
with the weak topology, such that for any $\bs\in\BN$, if
\[
\#\big\{i: \hat\lambda_i(\bs) \leq -2\nu''(q)^{\frac{1}{2}}+\epsilon/2 \big\} < N \delta ,
\]
then $\mu_{\bs}\notin A$ where $\mu_{\bs}$ denotes the empirical 
measure 
\[
\mu_{\bs}:=\frac{1}{N-1}\sum_{i=1}^{N-1}\delta_{\tilde\lambda_{i}(\bs)}\,.
\]

Hence, by the LDP for the empirical measure of GOE matrices \cite[Theorem 1.1]{BAG97},
for some $a=a(\epsilon,\delta)>0$, 
\eqref{eq:unifout1} holds with a bound of $e^{-N^{2}a}$ instead of $e^{-N\hat c}$, which is even stronger, for large $N$.

Since, for any $\bs\in\BN$, $\tilde\lambda_K(\bs)$  has the same distribution as the $K$-th smallest eigenvalue of a GOE matrix, 
by \cite[Theorem A.9]{A-BA-C}, it
satisfies a LDP at speed $N$ with rate function given by $KJ(t)$
where for $t> -2$, $J(t)=\infty$ and for $t\leq -2$,
\[
J(t)=\int_{t}^{-2}\sqrt{\frac{1}{4}s^{2}-1}ds.
\]
Thus, for any $\hat c>0$ and $\epsilon>0$, if
$K$ is large enough, then
\begin{equation*}
\sup_{\bs \in \BN}\P\Big\{ \#\big\{i:\tilde \lambda_i(\bs) \leq -2-\epsilon/\nu''(q)^{\frac{1}{2}} \big\} >K\Big\} \leq e^{-N\hat{c}},
\end{equation*}
from which \eqref{eq:unifout2} follows. This completes the proof of the lemma. \qed

\subsection{Proof of Theorem \ref{thm:LB}}
Suppose that $\bv_j$ and $\bs_{j/k}$ are as in Part \eqref{Item1} of the theorem. On the event from \eqref{eq:Cont}, which occurs with probability at least $1-e^{-NC}$ for appropriate $C>0$, by Taylor's theorem, for any $t\in[0,1/k]$ and $1\leq j\leq k-1$, 
\begin{equation*}
\frac{1}{N}H_{N}(\bs_{j/k+t})<\frac{1}{N}H_{N}(\bs_{j/k})-\nu''(j/k)^{\frac{1}{2}}t+\frac{\epsilon}{2k}+\frac{R}{6k^{3/2}},
\end{equation*}
and for $j=0$,
\begin{equation*}
\frac{1}{N}H_{N}(\bs_{t})<R/\sqrt{k}.
\end{equation*}
Part \eqref{Item1} therefore  easily follows with an appropriate choice of $\eta$, depending on the bounds above and the modulus of continuity of $\nu''(q)$ on $[0,1]$.

For Part \eqref{Item2}, consider the following optimization algorithm, written in pseudo-code. Define the coefficient $\tilde J_{i_{1},...,i_{p}}^{(p)}$ as the sum of $J_{i_{1}',...,i_{p}'}^{(p)}$ over all permutations $i_{1}',...,i_{p}'$ of $i_{1},...,i_{p}$.
\vspace{.15cm}
	
\begin{algorithm}[H]
		\caption*{Hessian based optimization}
		\SetAlgoLined
		\KwIn{All the disorder coefficients $J_{i_{1},...,i_{p}}^{(p)}$, parameters $\epsilon>0$, $k\geq1$}
		\KwOut{Sequence of directions $\bv_{j}$\ \ s.t. (w.h.p.)  $\frac1N H_N(\bs_{q})\leq -\EH(q) + \eta$ }
		Compute all the coefficients $\tilde J_{i_{1},...,i_{p}}^{(p)}$\;
		Initialize $\bv_0=(1,0,\ldots,0)$,  $\bs_{1/k}=(\sqrt{N/k},0,\ldots,0)$\;
		\For{$i= 1$ \KwTo $k-1$}{
			$\bs=\bs_{i/k}$, $q=\|\bs\|^2/N$\;
			Compute $\grad H_N(\bs)$ and $\Hess H_N(\bs)$\;
			Find   $\bv\perp \bs$, $\|\bv\|=1$, s.t. (w.h.p.): (1) $\bv^{T}\Hess H_{N}(\bs)\bv\leq -2\nu''(q)^\frac{1}{2}+\epsilon$,\newline 
			\phantom{Find $\bv\perp \bs$, $\|\bv\|=1$, s.t. (w.h.p.):}(2) $\bv^T\grad H_{N}(\bs)\leq0$\;
			$\bv_i = \bv$,\ \  $\bs_{i+1/k} =\bs + \sqrt{N/k}\bv$\;
			}
			\KwRet $(\bv_{0},\bv_{1},\ldots,\bv_{k-1})$
			\end{algorithm}
			
			\vspace{.15cm}Next, we explain how a direction $\bv$ as in the for-loop is found, and discuss the time complexity, assuming that $\deg(\nu)<\infty$ (see Footnote \ref{ft:degnu}). 
			
			Computing one element of $\nabla_E^2 H_N(\bs)$, 
			\[
			\frac{d}{dx_{i}}\frac{d}{dx_{j}}H_{N}(\bs)=\sum_{p=2}^{\infty}\frac{\gamma_{p}}{N^{\frac{p-1}{2}}}\sum_{i_1\leq \cdots\leq i_{p-2}}(1+w_i+\delta_{ij})(1+w_j)\tilde J_{i,j,i_{1},...,i_{p-2}}^{(p)}\sigma_{i_{1}}\cdots\sigma_{i_{p-2}},
			\] 
			takes $O(N^{\deg(\nu)-2})$ operations, where, abusing notation, $w_i=\#\{k:i_k=i\}$.
			From the Euclidean Hessian, $\Hess H_N(\bs)$ can be directly computed. The gradient can be treated similarly. It is therefore easy to see that the second line in the for-loop above takes $O(N^{\deg(\nu)})$ operations. 
			
			By Lemma \ref{lem:good_eigvals}, for small $\delta$ and $C$, 
			\begin{equation}
			\label{eq:subspace}
			\P\Big\{\,\forall\bs \in \BN:\: \,  \#\big\{i: \lambda_i(\bs) \leq -2\nu''(q)^{\frac{1}{2}}+\epsilon/2 \big\} \geq N \delta\Big\}\geq1-e^{-NC}	.
			\end{equation}
			By a simple variation of the proof of the same lemma, for large enough $L>0$ and small $C$,
			\begin{equation}
			\label{eq:absvalbd}
			\P\Big\{\,\forall\bs \in \BN,\, \forall i:\: \,  |\lambda_i(\bs)| < L\,  \Big\}\geq1-e^{-NC}	.
			\end{equation}
			
			For a uniformly chosen direction $\bu$ with $\bu\perp \bs$ and $\|\bu\|=1$, for small $C$, with probability at least $1-e^{-NC}	$, the projection of $\bu$ onto the span of the eigenvectors that correspond to the eigenvalues from \eqref{eq:subspace} is at least $\sqrt{\delta}/2$. From simple linear algebra, if this happens and the event from \eqref{eq:absvalbd} occurs, then for large enough integer $m=m(\delta,\epsilon,L)$, the power iteration 
			$$
			\bv = \frac{(\Hess H_N(\bs)-L\mathbf{I})^m\bu}{\|(\Hess H_N(\bs)-L\mathbf{I})^m\bu\|},
			$$
			where $\mathbf{I}$ is the identity matrix, satisfies  condition (1) from the algorithm. By flipping its sign if needed, we obtain a vector $\bv$ that also satisfies condition (2). Combining the above we have that the algorithm has time complexity $C'N^{\deg(\nu)}$, and it produces a sequence $\bv_{j}$ as required, with probability at least $1-e^{-NC}$, for small $C$. \qed

\subsection{Proof of Corollary \ref{cor:GSlb}}
Fix some $k\geq1$ and $\epsilon>0$. 
By Lemma \ref{lem:good_eigvals}, w.h.p. we can choose a set of directions $\bv_j$ as in Part \eqref{Item1} of Theorem \ref{thm:LB} one after the other. Namely, after choosing $\bv_0,\ldots,\bv_{j-1}$, which define $\bs_{j/k}$, there exists w.h.p. a direction $\bv_j$ orthogonal to $\bs_{j/k}$,  
which satisfies the second inequality of  \eqref{eq:vj}, and by flipping its sign if needed, also the first inequality of \eqref{eq:vj}.

Hence, by Part \eqref{Item1} of Theorem \ref{thm:LB}, for any $\eta>0$ and $q\in[0,1]$, we have that $-\Es(q)\leq -\EH(q)+\eta$.\qed

\subsection{Proof of Lemma \ref{lem:orthoGS}}

Let $k\geq1$ be some large integer number and let $a>0$ be some small error. For any model satisfying the Ghirlanda-Guerra identities, and any generic
mixed model in particular \cite[Section 3.7]{PanchenkoBook}, by
the ultrametricity property \cite{ultramet}, we have the following (see Remark 2.1 in \cite{PanchenkoBook}). 

If $R$ is a $k\times k$ overlap array whose diagonal elements are equal to $1$, 
its off-diagonal elements belong to the support of the Parisi measure, 
and it is ultrametric in the sense that $R_{ij}\geq\min\{R_{il},\,R_{jl}\}$, then 
the probability that the overlap array of $k$ points sampled independently from the Gibbs measure is equal to $R$, up to error $a$ uniform in the elements of the array,
is bounded away from $0$ for large $N$.

If a model does not satisfy
the Ghirlanda-Guerra identities, by approximating it by a sequence
of models that do satisfy the identities, we can obtain the same,
but with probability which only does not decay exponentially fast
in $N$ instead of being bounded away from $0$. 

This is explained in the proof of \cite[Lemma 24]{FElandscape} with
the overlap array being constant off the diagonal ($R_{ij}=q+\delta_{ij}(1-q)$).
However, the same proof works for arrays as above and we refer the
reader there for details. In particular, for general $\nu(x)$ assuming
that $q<q'\in S_{P}\cap(0,1)$, we obtain that for some sequences
$a_{N},\,c_{N}\to0$ and $k_{N}\to\infty$,
\begin{equation}
\lim_{N\to\infty}\P\left(\frac{1}{N}\log G_{N,\beta}^{\otimes\infty}\left\{ \forall i,j,i',j'\leq2k_{N}:\,\big|R(\bs_{i,j},\bs_{i',j'})-q_{i,j,i',j'}\big|<a_{N}\right\} >-c_{N}\right)=1,\label{eq:samples}
\end{equation}
where $G_{N,\beta}^{\otimes\infty}$ denotes the law of an infinite
array of independent samples $\{\bs_{i,j}\}_{i,j=1}^{\infty}$ from
the Gibbs measure $G_{N,\beta}$, and 
\[
q_{i,j,i',j'}=\begin{cases}
q & ,\text{ if }i\neq i',\\
q' & ,\text{ if }i=i',\,j\neq j',\\
1 & ,\text{ if }i=i',\,j=j'.
\end{cases}
\]

Next, we wish to average over `clusters' of the samples as in (\ref{eq:samples})
(which belong to $\SN(1)$) in order to obtain points in $\SN(q)$
and $\SN(q')$ which have overlaps as in (\ref{eq:OlapBd}), such
that spherical `bands' around them, i.e., sets of the form 
\begin{equation}
B(\bs,\rho)=\big\{\bs'\in\SN:\,|R(\bs,\bs')-R(\bs,\bs)|<\rho\big\},\label{eq:band-1}
\end{equation}
with small $\rho>0$, have Gibbs weight that does not decay exponentially
fast with $N$, and such that under the Gibbs measure conditional
on $B(\bs,\rho)$ the probability that many samples all have pair-wise
overlap roughly $R(\bs,\bs)$ is not exponentially small. This property
was shown in \cite[Proposition 10]{FElandscape} to imply that $H_{N}(\bs)/N$
is roughly equal to $-\Es(R(\bs,\bs))$, and this will therefore imply
the bounds of (\ref{eq:ValBd}).

More precisely, for $1\leq i\leq k_{N}$, we define 
\[
\hat{\bs}_{i}=\frac{1}{k_{N}}\sum_{j=1}^{k_{N}}\bs_{i,j}\text{\,\,\, and \,\,}\hat{\bs}_{0}=\frac{1}{k_{N}}\sum_{i=1}^{k_{N}}\hat{\bs}_{i}=\frac{1}{k_{N}^{2}}\sum_{i=1}^{k_{N}}\sum_{j=1}^{k_{N}}\bs_{i,j},
\]
and
\begin{equation}
\bs_{i}=\sqrt{Nq'}\frac{\hat{\bs}_{i}}{\|\hat{\bs}_{i}\|}\text{\,\,\,\,\,and\,\,\,\,}\bs=\sqrt{Nq}\frac{\hat{\bs}_{0}}{\|\hat{\bs}_{0}\|}.\label{eq:barbs}
\end{equation}
Note that assuming that the overlaps of $\bs_{i,j}$ satisfy the bounds
as in (\ref{eq:samples}), for $1\leq i\leq k_{N}$ and $1\leq i',\,j'\leq2k_{N}$,
\begin{equation}
\begin{aligned}R(\hat{\bs}_{i},\hat{\bs}_{0}),\,R(\hat{\bs}_{0},\hat{\bs}_{0}),\,R(\hat{\bs}_{0},\bs_{i',j'}) & \in[q-a_{N}',q+a_{N}'],\\
R(\hat{\bs}_{i},\hat{\bs}_{i'}),\,R(\hat{\bs}_{i},\bs_{i',j'}) & \in[q_{i,i'}-a_{N}',q_{i,i'}+a_{N}'],
\end{aligned}
\label{eq:overlap_relation}
\end{equation}
and
\begin{equation}
\begin{aligned}R(\bs_{i},\bs),\,R(\bs,\bs),\,R(\bs,\bs_{i',j'}) & \in[q-\rho_{N},q+\rho_{N}],\\
R(\bs_{i},\bs_{i'}),\,R(\bs_{i},\bs_{i',j'}) & \in[q_{i,i'}-\rho_{N},q_{i,i'}+\rho_{N}],
\end{aligned}
\label{eq:overlap_relation-1}
\end{equation}
where $a_{N}'=a_{N}+k_{N}^{-1}$ and, for large $N$, $\rho_{N}=a_{N}'+5\sqrt{a_{N}'}$
and where
\[
q_{i,i'}=\begin{cases}
q & ,\text{ if }i\neq i',\\
q' & ,\text{ if }i=i'.
\end{cases}
\]

Therefore, on the event in (\ref{eq:samples}), by conditioning on
the samples $\bs_{i,j}$ from $G_{N,\beta}^{\otimes\infty}$ with
$1\leq i,\,j\leq k_{N}$, we have that with high probability:
\begin{enumerate}
\item \label{enu:p1}The points $\bs_{i}$ and $\bs$
defined by (\ref{eq:barbs}) satisfy (\ref{eq:overlap_relation-1}).
\item \label{enu:p2}For the band around $\bs_i$ for any $1\leq i\leq k_{N}$, 
\begin{align*}
\frac{1}{N}\log G_{N,\beta}^{\otimes\infty}\Big\{ & \,\forall k_{N}<j,\,j'\leq2k_{N},\,j\neq j':\\
 & \bs_{i,j}\in B(\bs_{i},\rho_{N}),\,\big|R(\bs_{i,j},\bs_{i,j'})-q'\big|<\rho_{N}\Big\}>-c_{N}.
\end{align*}
\item \label{enu:p3}For  the band around $\bs$,
\begin{align*}
\frac{1}{N}\log G_{N,\beta}^{\otimes\infty}\Big\{ & \,\forall k_{N}<j,\,j'\leq2k_{N},\,j\neq j'\\
 & \bs_{j,j}\in B(\bs,\rho_{N}),\,\big|R(\bs_{j,j},\bs_{j',j'})-q\big|<\rho_{N}\Big\}>-c_{N}.
\end{align*}
\end{enumerate}
By \cite[Proposition 10]{FElandscape}, for some $t_{N}=o(1)$, on
an event whose probability tends to $1$ as $N\to\infty$, the points
$\bs_{i}\in\SN(q')$ in Point (\ref{enu:p2}) above also satisfy
\[
\left|\frac{H_{N}(\bs_{i})}{N}+\Es(q')\right|<t_{N},
\]
and the point $\bs\in\SN(q)$ in Point (\ref{enu:p3}) also
satisfies
\[
\left|\frac{H_{N}(\bs)}{N}+\Es(q)\right|<t_{N}.
\]

Combined with Point (\ref{enu:p1}), this means that with probability
going to $1$, the points $\bs$ and $\bs_{i}$
satisfy (\ref{eq:OlapBd}) and (\ref{eq:ValBd}), and the proof is completed.\qed

\subsection{Proof of Proposition \ref{prop:upper bound}}

We claim that we only need to show that 
\begin{align}
\limsup_{\epsilon\searrow0}\frac{\Es(q+\epsilon)-\Es(q)}{\epsilon} & \leq\nu''(q)^{1/2}.\label{eq:dEsUB}
\end{align}
Indeed, using the same argument as
in Theorem \ref{thm:LB} to construct a path starting from the minimizer
\[
\bs_{q}:=\argmin_{\bs\in\SNq}H_{N}(\bs),
\]
instead of the origin $\bs_{0}=0$, we conclude similarly to Corollary
\ref{cor:GSlb} that for small $\epsilon>0$,
\[
\Es(q+\epsilon)-\Es(q)\geq\int_{q}^{q+\epsilon}\nu''(t)^{1/2}dt.
\]
Therefore,
\begin{align*}
\liminf_{\epsilon\searrow0}\frac{\Es(q+\epsilon)-\Es(q)}{\epsilon} & \geq\nu''(q)^{1/2},
\end{align*}
and if we assume  (\ref{eq:dEsUB}) then (\ref{eq:rightder}) follows. From Lemma \ref{lem:Cont} it is easy to see that $\Es(q)$ is a continuous function of $q\in[0,1]$. Thus, if $[0,q]\subset S_{P}$, from (\ref{eq:rightder})
we have that $\Es(q)=\int_{0}^{q}\nu''(t)^{1/2}dt$.


Let $\epsilon>0$ be some arbitrary small number such that $q,\,q+\epsilon\in S_{P}$.
From now on, we will assume that the following events, which have
probability tending to $1$ as $N\to\infty$, occur. First, we assume
that the bounds on the derivatives as in  \eqref{eq:Cont} hold.
Second, relying on Lemma \ref{lem:orthoGS}, we assume for appropriate
$\delta_{N}\to0$ and $k_{N}\to\infty$, $\bs\in\SNq$ and $\bs_{1},\ldots,\bs_{k_{N}}\in\SN(q+\epsilon)$
are points satisfying the bounds of (\ref{eq:OlapBd}) and (\ref{eq:ValBd})
with $q'=q+\epsilon$.

We now fix $\epsilon$  and keep track
of the asymptotic behavior only in $N$. Namely, $o(1)$ will stand
for terms going to $0$ as $N\to\infty$ and $\epsilon$ is fixed,
with the rate in the bound being deterministic (on the event we have
restricted to). 

With 
\[
\bu_i:=\frac{\bs_{i}-\bs}{\|\bs_{i}-\bs\|},\,\,\,\text{where}\,\,\,\|\bs_{i}-\bs\|=\sqrt{N}\left(\sqrt{\epsilon}+o(1)\right),
\]
for all $i\leq k_{N}$,
\begin{equation}
\begin{aligned} & \Es(q+\epsilon)-\Es(q)\leq-\frac{1}{N}\left(H_{N}(\bs_{i})-H_{N}(\bs)\right)+o(1)\\
 & \leq-\sqrt{\frac{\epsilon}{N}}\grad H_{N}(\bs)\cdot\bu_i-\frac{\epsilon}{2}\bu_i^{T}\grad^{2}H_{N}(\bs)\bu_i+\frac{R}{6}\epsilon^{\frac{3}{2}}+o(1).
\end{aligned}
\label{eq:dEsBasicBd}
\end{equation}
Note that in the usual Taylor expansion we should have $\nabla_{E}$
and $\nabla_{E}^{2}$ instead of $\grad$ and $\Hess$ as above. But
since the projection of $\bu_i$ onto the direction of $\bs$ is small by (\ref{eq:OlapBd}), using $\grad$ and $\Hess$
instead results in an error which we absorbed into the $o(1)$ term.
Note that by the bound we assumed on the derivatives of $H_{N}(\bs)$,
\[
\left|\sqrt{\frac{1}{N}}\grad H_{N}(\bs)\cdot\bu_i\right|\leq R\left|\frac{\grad H_{N}(\bs)}{\|\grad H_{N}(\bs)\|}\cdot\bu_i\right|.
\]

Thus, to prove (\ref{eq:dEsUB}) and finish the proof, it will be
enough to show that for arbitrary fixed $\epsilon,\,c>0$, with probability
bounded away from $0$ uniformly in large $N$, for at least one of
the directions $\bu_i$ satisfying (\ref{eq:dEsBasicBd}),
we have that
\begin{equation}
R\left|\frac{\grad H_{N}(\bs)}{\|\grad H_{N}(\bs)\|}\cdot\bu_i\right|<c\sqrt{\epsilon},\label{eq:bd1}
\end{equation}
and 
\begin{equation}
\bu_i^{T}\grad^{2}H_{N}(\bs)\bu_i\geq-(2+c)\nu''(q)^{1/2}.\label{eq:bd2}
\end{equation}
Since $k_{N}\to\infty$, it will be enough to show that with probability
bounded away from $0$ uniformly in large $N$, (\ref{eq:bd1}) and
(\ref{eq:bd2}) do not occur for at most $K$ of the indices $i\leq k_{N}$,
for some $K$ independent of $N$.

Since the vectors $\bu_i$ are approximately orthogonal (\ref{eq:OlapBd}), by
Pythagoras' Theorem and simple linear algebra to account for the $\delta_{N}$
errors, it follows that for large enough $N$
there are at most 
\[
\frac{R^{2}}{c^{2}\epsilon}+1
\]
indices $i\leq k_{N}$ such that (\ref{eq:bd1}) does not occur.

If (\ref{eq:bd2}) does not occur for $K$ of the indices
$i\leq k_{N}$, then from (\ref{eq:OlapBd}) and the min-max theorem,
for large $N$,
\begin{equation}
\sqrt{\frac{N}{N-1}}\frac{\lambda_{K}(\bs)}{\nu''(q)^{1/2}}<-2-c/2.\label{eq:lambdaK}
\end{equation}
By Lemma \ref{lem:good_eigvals}, for large enough $K$ this occurs with probability going to $0$ as $N\to\infty$. Combining the above the proof is completed.
\qed

\subsection{Proof of Theorem \ref{thm:FEwoutParisi}}

From the lower bound of (\ref{eq:TAP}) (proved in \cite{FElandscape}) and Corollary \ref{cor:GSlb},
 (\ref{eq:FElb}) will follow if we show that
\[
F_{\beta}(q_{P})=\frac{1}{2}\beta^{2}\nu_{q_{P}}(1),
\]
where if $\beta\leq\nu''(0)^{-\frac{1}{2}}$ and $S_{P}=\{0\}$ we
define $q_{P}=0$, and otherwise $q_{P}$ is the unique solution of
$\nu''(q)^{-\frac{1}{2}}=\beta(1-q)$ (see Proposition \ref{prop:1}).

This will follow from Proposition \ref{prop:RSFE} if we can show
that
\[
g(s):=\beta^{2}\nu_{q_{P}}(s)+\log(1-s)+s<0,\quad \forall s\in(0,1),
\]
where we recall that
\begin{equation}
\nu_{q}(s)=\nu(q+(1-q)s)-\nu(q)-(1-q)\nu'(q)s\label{eq:nuq}
\end{equation}
(and where $\nu_{q}(s)=\nu(s)$ coincides with the original mixture when
$q=0$).

Since $g(0)=g'(0)=0$, it will be enough to show that 
\[
g''(s)=\beta^{2}(1-q_{P})^{2}\nu''(q_{P}+(1-q_{P})s)-(1-s)^{-2}
\]
is strictly negative for any $s\in(0,1)$. Note that for such $s$
the sign of $g''(s)$ is the same as that of
\[
h(s):=\beta(1-q_{P})(1-s)-\nu''(q_{P}+(1-q_{P})s)^{-\frac{1}{2}}.
\]
By the assumption on $\nu''(s)^{-1/2}$, $h(s)$ is a convex function
with $h(1)<0$. For $\beta>\nu''(0)^{-\frac{1}{2}}$, since $\beta(1-q_{P})\nu''(q_{P})^{\frac{1}{2}}=1$,
$h(0)=0$. And for $\beta\leq\nu''(0)^{-\frac{1}{2}}$, $q_{P}=0$
and thus $h(0)\leq0$. From convexity we have that $h(s)<0$, and
therefore $g''(s)<0$, if $s\in(0,1)$ and we conclude (\ref{eq:FElb}).

When $q_{P}=0$ and $x_{P}\equiv1$, clearly the right-hand side of
(\ref{eq:FElb}) coincides with $\mathcal{P}(x_{P})$. Assume that
$\beta>\nu''(0)^{-\frac{1}{2}}$ and therefore $q_{P}>0$. By substituting
(\ref{eq:eta}) we obtain that
\begin{equation}
\begin{aligned}\mathcal{P}(x_{P}) & =\frac{1}{2}\Big(\beta\int_{0}^{q_{P}}\eta(q)\nu'(q)dq\\
 & +\beta^{2}(\nu(1)-\nu(q_{P}))+\beta\int_{0}^{q_{P}}\nu''(q)^{1/2}dq+\log(1-q_{P})\Big),
\end{aligned}
\label{eq:PxP}
\end{equation}
where we used the fact that for $q\leq q_{P}$,
\[
\hat{x}_{P}(q):=\int_{q}^{1}x(s)ds=\frac{1}{\beta}(\nu''(q)^{-1/2}-\nu''(q_{P})^{-1/2})+(1-q_{P})=\frac{1}{\beta}\nu''(q)^{-1/2}.
\]

By integration by parts,
\begin{equation}
\begin{aligned}\int_{0}^{q_{P}}\eta(q)\nu'(q)dq & =\int_{0}^{q_{P}}\nu''(q)^{1/2}dq-\frac{\nu'(q_{P})}{\nu''(q_{P})^{1/2}}\\
 & =\int_{0}^{q_{P}}\nu''(q)^{1/2}dq-\beta(1-q_{P})\nu'(q_{P}).
\end{aligned}
\label{eq:intbypts}
\end{equation}
From (\ref{eq:nuq}), (\ref{eq:PxP}) and (\ref{eq:intbypts}) we
obtain that $\mathcal{P}(x_{P})$ is equal to the right-hand side
of (\ref{eq:FElb}).\qed

\subsection*{Acknowledgment}

The author was supported by the Simons Foundation. This project has received funding from the European Research Council (ERC) under the European Union's Horizon 2020 research and innovation programme (grant agreement No. 692452).

\bibliographystyle{plain}
\bibliography{master}

\end{document}